\documentclass[10pt,a4paper,oneside]{article}
\usepackage[latin1]{inputenc}
\usepackage{amsmath}
\usepackage{amsfonts}
\usepackage{amssymb}
\usepackage{latexsym}

\newtheorem{theorem}{Theorem}[section]

\newtheorem{lemma}[theorem]{Lema}

\def \fimprova {\blacksquare}

\begin{document}

\title{Classifying Finitely Generated Indecomposable RA Loops}
\author{Mariana G. Cornelissen \\ \small\emph{{Universidade Federal de São João Del Rei, Brazil}}
\\ \\ C.Polcino Milies \\ \small\emph{{Universidade de São Paulo, Brazil}}}
\date{}

\maketitle

\begin{abstract}
In \cite{casofinito}, E. Jespers, G. Leal and C.
Polcino Milies classified all finite ring alternative loops (RA loops for short) which are not direct products of proper subloops. In this paper we extend this result to finitely generated RA loops and provide an explicit description of all such loops.
\end{abstract}

\section{Basic Definitions}

A loop is a pair $(L,.)$ where $L$ is a nonempty set, $(a,b)
\mapsto a.b$ is a closed binary operation on $L$ which has a two
sided identity element $1$ and with the property that the equation
$a.b=c$ determines a unique element $b \in L$ for given $a,c \in
L$ and a unique element $a \in L$ for given $b,c \in L$. A subloop of a loop $(L,.)$ is a subset of $L$ which, under the binary operation, is also a loop.


Given a commutative and associative ring $R$ with unity and a loop
$L$, we can construct the loop algebra $RL$ of $L$ over $R$ as the
free $R$-module with basis $L$ in which the multiplication is
defined by extending that of $L$ via the distributive laws.

We remember that a ring $R$ (not necessarily associative) is called alternative if it
satisfies the identities
$$ [x,x,y]=0 \mbox{ and } [y,x,x]=0$$
for all $x,y \in R$ where $[a,b,c]=(ab)c - a(bc)$ is the \emph{associator} of $a,b,c$.

A \textbf{R}ing \textbf{A}lternative Loop (RA loop) is a loop whose loop ring over some
commutative, associative ring with unity and of characteristic
different from 2, is alternative but not associative.

In this paper we provide an explicit description of all
finitely generated indecomposable RA loops, i.e.,
finitely generated RA loops which are not direct products of proper subloops.

\section{Some Results}

We list below some results about loops and groups that are used in
this work. These proofs can be found in \cite{livro}, \cite{artigogrupos}.



Suppose that $G$ is a non abelian group with center $Z(G)$, $g_o$
is an element of $Z(G)$ and that $g \mapsto g^*$ is an involution
in $G$ such that ${g_{o}}^* =g_o$ and $gg^* \in Z(G)$ for every $g
\in G$. Let $L = G \cup Gu$ where $u$ is an indeterminate and
extend the binary operation from $G$ to $L$ by the rules:
\[ g(hu)=(hg)u \]
\[ (gu)h=(gh^*)u \]
\[ (gu)(hu)=g_oh^*g \]
for $g,h \in L$. It can be shown that $L$ is a non-associative
loop where each element has a unique two-sided inverse. This loop
is denoted by $M(G,*,g_o)$. The theorem below shows that the loops
obtained in this way are, in fact, all the RA loops.

\begin{theorem}[\cite{livro}, Theorem IV.3.1] \label{raloops}
Let $L$ be a RA loop. Then there exists a group $G \subset L$ and an
element $u \in L$ such that $L= G \cup Gu $, $L'=G'= \{1,s\} \subseteq Z(G) = Z(L)$, where $L'$ is the commutator-associator subloop of $L$, $L / Z(L) \simeq C_2 \times C_2 \times
C_2$, where $C_2$ denotes a cyclic group of order 2,
and consequently $G / Z(G) \simeq C_2 \times C_2$. \\
Furthermore, the map $*: L \rightarrow L$ given by

\begin{eqnarray}\label{involution}
 g^* = \left\{%
\begin{array}{ll}
$g$, & \hbox{if $g \in Z(G)$;} \\
$sg$, & \hbox{if $g \notin Z(G)$.} \\
\end{array}%
\right.
\end{eqnarray}
is an involution of $L$, and, setting $u^2 =g_o$, we have that
$g_o \in Z(G)$ and $L = M(G,*,g_o)$.

Conversely, for any nonabelian group $G$ such that $G / Z(G) \simeq C_2
\times C_2$ and involution *, the loop $M(G,*,g_o)$ is a RA loop
for any $g_o \in Z(G)$, where * is given by (\ref{involution}).
\end{theorem}

According to the theorem above, in order to classify the finitely generated RA loops, we have to study the finitely generated groups $G$ such that $G/Z(G) \simeq C_2 \times C_2$. In \cite{artigogrupos}, M.Cornelissen and C.P.Milies classified all finitely generated groups $G$ such that $G/Z(G) \simeq C_p \times C_p$ where $C_p$ denotes a cyclic group of prime order $p$. We quote below two important results of this work.

\begin{theorem}[\cite{artigogrupos}, Theorem 2.1] \label{decomposition}
A finitely generated group $G$ is such that $G/Z(G)\simeq C_p \times
C_p$, where $C_p$ denotes a cyclic group of prime order $p$, if and only if it can be written in the form $G = D \times A$
where $A$ is a finitely generated abelian group and $D$ is an  indecomposable group such that $D = \langle x,y,Z(D)\rangle$ where
$Z(D)$ is of the form
 $Z(D) = \langle t_1\rangle \times \langle z_2\rangle \times \langle z_3\rangle$, with:

 $(i)$  $o(t_1)= p^{m_1}, m_1 \geq 1$ and $s = [x,y]=x^{-1}y^{-1}xy = t_1^{p^{m_1-1}}$,

 $(ii)$ either $o(z_i)=
 p^{m_i}$ with $m_i \geq 0$ or $o(z_i) = \infty$, for $i=2,3$,

 $(iii)$
$x^p \in \langle t_1\rangle
\times \langle z_2\rangle$ and $y^p \in \langle t_1\rangle \times \langle z_2\rangle \times \langle z_3\rangle$.
\end{theorem}


\begin{theorem}[\cite{artigogrupos}, Theorem 3.8]\label{groupclassification}
     Let $G$ be a finitely generated indecomposable group such that $G/Z(G) \simeq C_p \times C_p$
     where $C_p$ denotes a cyclic group of prime order $p$. Then $G$ is of the form $G = \langle x,y, Z(G) \rangle$
     with $x^p, y^p \in Z(G)$ and belongs to one of the nine types of non isomorphic groups listed in the table below.
     In each case we have that $o(t_i) = p^{m_i}, i=1,2,3$ and $o(u_j) = \infty, j = 1,2$.


\begin{center}
\begin{table}[htb]
\begin{center}
\caption{Classification Group Table}
\medskip
\begin{tabular}{||c||c||c||c||}
  \hline
  $G$ & $Z(G)$ & $x^p$ & $y^p$ \\
  \hline
  1 & $\langle t_1 \rangle$ & 1 & 1 \\
  2 & $\langle t_1 \rangle$ & $t_1$ & $t_1$ \\
  3 & $\langle t_1 \rangle \times \langle t_2 \rangle$ & 1 & $t_2$ \\
  4 & $\langle t_1 \rangle \times \langle t_2 \rangle$ & $t_1$ & $t_2$ \\
  5 & $\langle t_1 \rangle \times \langle u_1 \rangle$ & 1 & $u_1$ \\
  6 & $\langle t_1 \rangle \times \langle u_1 \rangle$ & $t_1$ & $u_1$ \\
  7 & $\langle t_1 \rangle \times \langle t_2 \rangle \times \langle t_3 \rangle$ & $t_2$& $t_3$ \\
  8 & $\langle t_1 \rangle \times \langle t_2 \rangle \times \langle u_1 \rangle$  & $t_2$ & $u_1$ \\
  9 & $\langle t_1 \rangle \times \langle u_1 \rangle \times \langle u_2 \rangle$ & $u_1$ & $u_2$ \\
  \hline
\end{tabular}
\end{center}
\end{table}
\end{center}
\end{theorem}

Hence, in this work, we will use the above classification with $p=2$ to describe all the finitely generated indecomposable loops of the type $L=M(G,*,g_o)$.


\section{Classifying Finitely Generated Indecomposable RA Loops}

With the results of the previous section and a few other remarks
we shall show how to construct all the indecomposable finitely
generated RA loops.

\begin{lemma}[\cite{livro}, Proposition V.1.6]\label{lema}
Let $L = M(G,*,g_o)$ an RA loop and $A$ an abelian group. Then, $(g,a)^* = (g^*,a)$ defines an involution in $G \times A$ and $M(G \times A, *,(g_o,1))$ is an RA loop isomorphic to $M(G,*,g_o) \times A$. Conversely, if $A$ is an abelian group such that $M(G \times A, *, (g_o,1))$ is an RA loop for some non abelian group $G$ and some $g_o \in Z(G)$, then $*$ restricts to $G$ is an involution of $G$ and $M(G \times A, *,(g_o,1))\simeq M(G,*,g_o)\times A$.
\end{lemma}

We know by [\cite{livro}, Theorem IV.2.1] that squares are central in RA loops. Using the same techniques as in  [\cite{artigogrupos}, Lemma 2.2], it can be shown the following lemma that will be used frequently in what follows.

\begin{lemma}\label{expoente}
Let $v$ be an element of the indecomposable RA loop $L = M(G,*, g_o)$. Then $v$ can be chosen in such a way that if $w^{\alpha}$, $\alpha \in \mathbb{N}$, is a factor of $v^2 \in Z(L)$, we can write $\alpha = 0$ if $\alpha$ is even and $\alpha =1$, if $\alpha$ is odd.
\end{lemma}

%



Now we observe that if $L=M(G,*,g_o)$ is an indecomposable RA
loop, then $G$ is close to being itself indecomposable.

\begin{theorem}Let $L=M(G,*,g_o)$ be a finitely generated
indecomposable RA loop. Then $G = D \times H$ where $D$ is a
finitely generated indecomposable group and $H$ is a cyclic group. If $H = \langle h \rangle$ is non trivial,
then $g_o = dh$ with $d \in Z(D)$ and $h \in H$.
\end{theorem}

\emph{Proof:} Let $G = G \cup Gu$, $u \notin G$, $u^2=g_o \in
Z(G)$. By Theorem \ref{decomposition} we have that $G = D \times
A$ with $D = \langle x,y,Z(D) \rangle$ indecomposable, $rank[Z(D)] \leq 3$ and $A$ a finitely generated abelian group. If $A$ is
trivial, we have nothing to prove. Suppose $A \not= \{1\}$ and
write $g_o = da$, $d \in Z(D), a\in A$. We claim that $a \not= 1$. Otherwise, $g_o = d \in Z(D)$ and, by Lemma \ref{lema}, we can write $L = M(D \times A, *, (g_o,1)) \simeq M(G,*,g_o) \times A$ which contradicts the non decomposability of $L$.  We can
write $A = B \times C \times F$ where $B$ is a 2-group, $|C|$ is odd and $F$ is a free group of finite rank.

We can suppose that $a \in B \times F$. In fact, if $a = a'c$ with $a' \in B \times F$ and $c \in C$ such that $c^n =1$ then we can change $u$ to $u'= \gamma u$ where $\gamma \in C$ and such that $\gamma ^2 = c^{n-1}$. (Observe that there exists such $\gamma$ since the map $x \mapsto x^2$ is an automorphism of $C$). Hence, $(u')^2 = (\gamma u)^2 = \gamma ^2 u^2 = c^{n-1}g_o = c^{n-1}da = c^{n-1}da'c = da' \in D \times B \times F$. (Note that $L = G \cup Gu = G \cup Gu'$).

Write $A = \langle t_1 \rangle \times  \langle t_2 \rangle \times \ldots \times
 \langle t_k \rangle \times \langle x_1 \rangle \times \ldots \times \langle x_l \rangle \times C$ with
$o(t_i)= 2^{m_i}$, $o(x_j)= \infty$ and $g_o = dt_1^{a_1}\ldots t_m^{a_m} x_1^{b_1} \ldots x_n^{b_n}$ with $m \leq k, n \leq l$ .

Using Lemma \ref{expoente} for $u^2 = g_o$ and remembering that $Z(G) = Z(L)$ and $L = G \cup Gu = G \cup G(\alpha u)$ with $\alpha \in Z(G)$, we can assume that $a_i, b_j \in \{0,1\}$ for $i = 1, \ldots ,m$ and $j = 1, \ldots n$, i.e, $g_o = dt_1\ldots t_m x_1 \ldots x_n$, $m \leq k, n \leq l$ and reordering, if necessary, we can suppose that
$o(t_1) \geq o(t_2) \geq \ldots \geq o(t_m)$.

Let $H = \langle t_1\ldots t_m x_1 \ldots x_n \rangle$ be a cyclic group. If $x_1 \ldots x_n =1$ then
$H = \langle t_1 \ldots t_m \rangle, g_o = dt_1 \ldots t_m$ and $A = \langle t_1\ldots
t_m \rangle \times \langle t_2\rangle \times \ldots \times \langle t_m \rangle \times \ldots \times
\langle t_k \rangle \times \langle x_1 \rangle \times \ldots \times \langle x_l \rangle \times C$. If
$x_1\ldots x_n \not= 1$ we can then write $A =  \langle t_1\ldots
t_mx_1\ldots x_n \rangle \times  \langle x_2 \rangle \times \ldots \times \langle x_n \rangle \times
\ldots \times \langle x_l \rangle \times \langle t_1 \rangle \times \ldots \times \langle t_k \rangle \times
C$. In both cases, we have that $A = H \times K$ with $H$ a cyclic group and $K$ an abelian group. Hence $ G = D
\times H \times K$ with $g_o \in D \times H$ and $K$ abelian.
Using again Lemma \ref{lema} we have that $L = M(D \times H, *,
g_o) \times K$. Indecomposable property of L implies that $K =
\{1\}$, so $ G = D \times H$ with $H$ cyclic, $H = \langle h \rangle$ where $g_o
= dh$. $\hfill{\fimprova}$

\medskip

Using the previous results, we can write the indecomposable RA loop $L$ in the form
$L = M(G,*,g_o) = G \cup Gu$ with:

$(i)$ $u^2 = g_o \in Z(G)$

$(ii)$ $G = D \times H = \langle x,y,Z(D) \rangle \times H$

$(iii)$ $Z(L) = Z(G) = Z(D) \times H = \langle t_1 \rangle \times \langle z_2 \rangle \times \langle z_3 \rangle \times \langle h \rangle$ where $o(t_1) = 2^{m_1}, m_1 \geq 1$, $o(z_i) = 2^{m_i}, m_i \geq 0$ or $o(z_i) = \infty$ for $i=2,3$ and $o(h) = 2^k, k \geq 0$ or $o(h) = \infty$

$(iv)$ $x^2, y^2 \in Z(D)$.

\medskip

Note that if $z$ is a non trivial factor of $x^2$ or $y^2$ then we can choose $u$ such that $u^2$ does not have
the factor $z$ in its decomposition, except in the case where $z = t_1$ and $m_1 =1$, i.e, ${t_1}^2=1$. In fact, if $u^2$ and $x^2$ has $z$ in its decomposition, set $ u' =xu$; then $L = G \cup Gu = G \cup Gu'$ and $u'^2 = sx^2u^2$ where $s$ is the unique non trivial commutator-associator of $L$. Since
$s = t_1^{2^{m_1-1}}$ the exponent of $z$ in $u'^2$ is equal two
(*unless $i=1$ and $m_1=1$). Hence, by Lemma \ref{expoente}, we can therefore make a new choice of
$u'$ such that $u'^2$ does not have $z$. A similar statement is true for $y^2$.

\medskip

We are now ready to enumerate all possible indecomposable RA
loops.
Using Theorem \ref{groupclassification} with $p=2$, it follows that exist 9 types of groups $D$ as above and for each one, we have six possibilities of indecomposable RA loops, as we can see in the tables below.

\begin{quote}
\begin{center}
Table 1: RA Loops when $D = G$ of type 1
\end{center}
\end{quote}

\begin{center}
\begin{tabular}{||c||c||c||c||c||c||}
  \hline
  Loop $L$ & $Z(D)$ & $x^2$ & $y^2$ & $G$ & $u^2=g_{o}$ \\
  \hline
  $L_1$ & $ \langle t_1 \rangle$ & 1 & 1 & $G_1$ & 1 \\
  $L_2$ & $\langle t_1 \rangle$ & 1 & 1 & $G_1$ & $t_1$ \\
  $L_3$ & $\langle t_1 \rangle$ & 1 & 1 & $G_1 \times \langle t \rangle$ & $t$ \\
  $L_4$ & $\langle t_1 \rangle$ & 1 & 1 & $G_1 \times \langle t \rangle$ & $t_1t$ \\
  $L_5$ & $\langle t_1 \langle$ & 1 & 1 & $G_1 \times \langle w \rangle$ & $w$ \\
  $L_6$ & $\langle t_1 \rangle$ & 1 & 1 & $G_1 \times \langle w \rangle$ & $t_1w$ \\
  \hline
\end{tabular}
\end{center}

\begin{quote}
\begin{center}
Table 2: RA Loops when $D = G$ of type 2
\end{center}
\end{quote}

\begin{center}
\begin{tabular}{||c||c||c||c||c||c||}
  \hline
  Loop $L$ & $Z(D)$ & $x^2$ & $y^2$ & $G$ & $u^2=g_{o}$ \\
  \hline
  $L_7$ & $\langle t_1 \rangle$ & $t_1$ & $t_1$ & $G_2$ & 1 \\
  $L_8^*$ & $\langle t_1 \rangle$ & $t_1$ & $t_1$ & $G_2$ & $t_1$ \\
  $L_9$ & $\langle t_1 \rangle$ & $t_1$ & $t_1$ & $G_2 \times \langle t \rangle$ & $t$ \\
  $L_{10}^*$ & $\langle t_1 \rangle$ & $t_1$ & $t_1$ & $G_2 \times \langle t \rangle$ & $t_1t$ \\
  $L_{11}$ & $ \langle t_1 \rangle$ & $t_1$ & $t_1$ & $G_2 \times \langle w \rangle$ & $w$ \\
  $L_{12}^*$ & $\langle t_1 \rangle$ & $t_1$ & $t_1$ & $G_2 \times \langle w \rangle$ & $t_1w$ \\
 \hline
  \end{tabular}
\end{center}

\begin{quote}
\begin{center}
Table 3: RA Loops when $D = G$ of type 3
\end{center}
\end{quote}

\begin{center}
\begin{tabular}{||c||c||c||c||c||c||}
  \hline
  Loop $L$ & $Z(D)$ & $x^2$ & $y^2$ & $G$ & $u^2=g_{o}$ \\
  \hline
  $L_{13}$ & $\langle t_1 \rangle \times \langle t_2 \rangle$ & 1 & $t_2$ & $G_3$ & 1 \\
  $L_{14}$ & $\langle t_1 \rangle \times \langle t_2 \rangle$ & 1 & $t_2$ & $G_3$ & $t_1$ \\
  $L_{15}$ & $\langle t_1\rangle \times  \langle t_2 \rangle$ & 1 & $t_2$ & $G_3 \times \langle t \rangle$ & $t$ \\
  $L_{16}$ & $\langle t_1 \rangle \times \langle t_2 \rangle$ & 1 & $t_2$ & $G_3 \times \langle t \rangle$ & $t_1t$ \\
  $L_{17}$ & $\langle t_1 \rangle \times  \langle t_2 \rangle$ & 1 & $t_2$ & $G_3 \times \langle w \rangle$ & $w$ \\
  $L_{18}$ & $\langle t_1\rangle \times \langle t_2 \rangle$ & 1 & $t_2$ & $G_3 \times \langle w \rangle$ & $t_1w$ \\
\hline
\end{tabular}
\end{center}

\begin{quote}
\begin{center}
Table 4: RA Loops when $D = G$ of type 4
\end{center}
\end{quote}

\begin{center}
\begin{tabular}{||c||c||c||c||c||c||}
  \hline
  Loop $L$ & $Z(D)$ & $x^2$ & $y^2$ & $G$ & $u^2=g_{o}$ \\
  \hline
  $L_{19}$ & $\langle t_1 \rangle \times \langle t_2 \rangle$ & $t_1$ & $t_2$ & $G_4$ & 1 \\
  $L_{20}^*$ & $\langle t_1 \rangle \times \langle t_2 \rangle$ & $t_1$ & $t_2$ & $G_4$ & $t_1$ \\
  $L_{21}$ & $\langle t_1 \rangle \times \langle t_2 \rangle$ & $t_1$ & $t_2$ & $G_4 \times \langle t \rangle$ & $t$ \\
  $L_{22}^*$ & $ \langle t_1 \rangle \times \langle t_2 \rangle$ & $t_1$ & $t_2$ & $G_4 \times \langle t \rangle$ & $t_1t$ \\
  $L_{23}$ & $\langle t_1 \rangle \times \langle t_2 \rangle$ & $t_1$ & $t_2$ & $G_4 \times \langle w \rangle$ & $w$ \\
  $L_{24}^*$ & $\langle t_1 \rangle \times \langle t_2 \rangle$ & $t_1$ & $t_2$ & $G_4 \times \langle w \rangle$ & $t_1w$ \\
\hline
\end{tabular}
\end{center}

\begin{quote}
\begin{center}
Table 5: RA Loops when $D = G$ of type 5
\end{center}
\end{quote}

\begin{center}
\begin{tabular}{||c||c||c||c||c||c||}
  \hline
  Loop $L$ & $Z(D)$ & $x^2$ & $y^2$ & $G$ & $u^2=g_{o}$ \\
  \hline
  $L_{25}$ & $\langle t_1 \rangle \times \langle u_1 \rangle$ & 1 & $u_1$ & $G_5$ & 1 \\
  $L_{26}$ & $\langle t_1 \rangle \times \langle u_1 \rangle$ & 1 & $u_1$ & $G_5$ & $t_1$ \\
  $L_{27}$ & $\langle t_1 \rangle \times \langle u_1 \rangle$ & 1 & $u_1$ & $G_5 \times \langle t \rangle$ & $t$ \\
  $L_{28}$ & $\langle t_1 \rangle \times \langle u_1 \rangle$ & 1 & $u_1$ & $G_5 \times \langle t \rangle$ & $t_1t$ \\
  $L_{29}$ & $\langle t_1 \rangle \times \langle u_1 \rangle$ & 1 & $u_1$ & $G_5 \times \langle w \rangle$ & $w$ \\
  $L_{30}$ & $\langle t_1 \rangle \times \langle u_1 \rangle$ & 1 & $u_1$ & $G_5 \times \langle w \rangle$ & $t_1w$ \\
\hline
\end{tabular}
\end{center}

\begin{quote}
\begin{center}
Table 6: RA Loops when $D = G$ of type 6
\end{center}
\end{quote}

\begin{center}
\begin{tabular}{||c||c||c||c||c||c||}
  \hline
  Loop $L$ & $Z(D)$ & $x^2$ & $y^2$ & $G$ & $u^2=g_{o}$ \\
  \hline
  $L_{31}$ & $\langle t_1 \rangle \times \langle u_1 \rangle$ & $t_1$ & $u_1$ & $G_6$ & 1 \\
  $L_{32}^*$ & $\langle t_1 \rangle \times \langle u_1 \rangle$ & $t_1$ & $u_1$ & $G_6$ & $t_1$ \\
  $L_{33}$ & $\langle t_1 \rangle \times \langle u_1 \rangle$ & $t_1$ & $u_1$ & $G_6 \times \langle t \rangle$ & $t$ \\
  $L_{34}^*$ & $\langle t_1 \rangle \times \langle u_1 \rangle$ & $t_1$ & $u_1$ & $G_6 \times \langle t \rangle$ & $t_1t$ \\
  $L_{35}$ & $\langle t_1 \rangle \times \langle u_1 \rangle$ & $t_1$ & $u_1$ & $G_6 \times \langle w \rangle$ & $w$ \\
  $L_{36}^*$ & $\langle t_1 \rangle \times \langle u_1 \rangle$ & $t_1$ & $u_1$ & $G_6 \times \langle w \rangle$ & $t_1w$ \\
  \hline
  \end{tabular}
\end{center}

\begin{quote}
\begin{center}
Table 7: RA Loops when $D = G$ of type 7
\end{center}
\end{quote}

\begin{center}
\begin{tabular}{||c||c||c||c||c||c||}
  \hline
  Loop $L$ & $Z(D)$ & $x^2$ & $y^2$ & $G$ & $u^2=g_{o}$ \\
  \hline
  $L_{37}$ & $\langle t_1 \rangle \times \langle t_2 \rangle \times \langle t_3 \rangle$ & $t_2$ & $t_3$ & $G_7$ & 1 \\
  $L_{38}$ & $\langle t_1 \rangle \times \langle t_2 \rangle \times \langle t_3 \rangle$ & $t_2$ & $t_3$ & $G_7$ & $t_1$ \\
  $L_{39}$ & $\langle t_1 \rangle \times \langle t_2 \rangle \times \langle t_3 \rangle$ & $t_2$ & $t_3$ & $G_7 \times \langle t \rangle$ & $t$ \\
  $L_{40}$ & $\langle t_1 \rangle  \times \langle t_2 \rangle  \times \langle t_3 \rangle$ & $t_2$ & $t_3$ & $G_7 \times \langle t \rangle$ & $t_1t$ \\
  $L_{41}$ & $\langle t_1 \rangle \times \langle t_2 \rangle \times \langle t_3 \rangle$ & $t_2$ & $t_3$ & $G_7 \times \langle w \rangle$ & $w$ \\
  $L_{42}$ & $\langle t_1 \rangle \times \langle t_2 \rangle \times \langle t_3 \rangle$ & $t_2$ & $t_3$ & $G_7 \times \langle w \rangle$ & $t_1w$ \\
\hline
\end{tabular}
\end{center}

\begin{quote}
\begin{center}
Table 8: RA Loops when $D = G$ of type 8
\end{center}
\end{quote}

\begin{center}
\begin{tabular}{||c||c||c||c||c||c||}
  \hline
  Loop $L$ & $Z(D)$ & $x^2$ & $y^2$ & $G$ & $u^2=g_{o}$ \\
  \hline
  $L_{43}$ & $\langle t_1 \rangle \times \langle t_2 \rangle \times \langle u_1 \rangle$ & $t_2$ & $u_1$ & $G_8$ & 1 \\
  $L_{44}$ & $\langle t_1 \rangle \times \langle t_2 \rangle \times \langle u_1 \rangle$ & $t_2$ & $u_1$ & $G_8$ & $t_1$ \\
  $L_{45}$ & $\langle t_1 \rangle \times \langle t_2 \rangle \times \langle u_1 \rangle$ & $t_2$ & $u_1$ & $G_8 \times \langle t \rangle$ & $t$ \\
  $L_{46}$ & $\langle t_1 \rangle \times \langle t_2 \rangle \times \langle u_1 \rangle$ & $t_2$ & $u_1$ & $G_8 \times \langle t \rangle$ & $t_1t$ \\
  $L_{47}$ & $\langle t_1 \rangle \times \langle t_2 \rangle \times \langle u_1 \rangle$ & $t_2$ & $u_1$ & $G_8 \times \langle w \rangle$ & $w$ \\
  $L_{48}$ & $\langle t_1 \rangle \times \langle t_2 \rangle \times \langle u_1 \rangle$ & $t_2$ & $u_1$ & $G_8 \times \langle w \rangle$ & $t_1w$ \\
\hline
\end{tabular}
\end{center}

\begin{quote}
\begin{center}
Table 9: RA Loops when $D = G$ of type 9
\end{center}
\end{quote}

\begin{center}
\begin{tabular}{||c||c||c||c||c||c||}
  \hline
  Loop $L$ & $Z(D)$ & $x^2$ & $y^2$ & $G$ & $u^2=g_{o}$ \\
  \hline
  $L_{49}$ & $\langle t_1 \rangle \times \langle u_1 \rangle \times \langle u_2 \rangle$ & $u_1$ & $u_2$ & $G_9$ & 1 \\
  $L_{50}$ & $\langle t_1 \rangle \times \langle u_1 \rangle \times \langle u_2 \rangle$ & $u_1$ & $u_2$ & $G_9$ & $t_1$ \\
  $L_{51}$ & $\langle t_1 \rangle \times \langle u_1 \rangle \times \langle u_2 \rangle$ & $u_1$ & $u_2$ & $G_9 \times \langle t \rangle$ & $t$ \\
  $L_{52}$ & $\langle t_1 \rangle \times \langle u_1 \rangle \times \langle u_2 \rangle$ & $u_1$ & $u_2$ & $G_9 \times \langle t \rangle$ & $t_1t$ \\
  $L_{53}$ & $\langle t_1 \rangle \times \langle u_1 \rangle \times \langle u_2 \rangle$ & $u_1$ & $u_2$ & $G_9 \times \langle w \rangle$ & $w$ \\
  $L_{54}$ & $\langle t_1 \rangle \times \langle u_1 \rangle \times \langle u_2 \rangle$ & $u_1$ & $u_2$ & $G_9 \times \langle w \rangle$ & $t_1w$ \\
  \hline
\end{tabular}
\end{center}

In each of the tables above, we are supposing that $o(t_i)=2^{m_i}, m_i \geq 1$, $o(u_j) = \infty, j = 1,2$, $o(t) = 2^k, k \geq 1$ and $o(w) = \infty$. The lines denote by $*$ are there from the exceptional case when $m_1 =1$, $s = t_1$ and $t_1^2 = 1$.

\section{Classification Theorem}

In this section, we will see that many of the RA loops listed in the tables of the previous section are isomorphic. The next theorem shows that, in fact, we have only 16 different types of finitely generated RA loops.

\begin{theorem} \label{classification}
Any finitely generated indecomposable RA loop is
in one of the sixteen types of loops listed in the table below.

\begin{quote}
\begin{center}
Table 10: Classification Loop Table
\end{center}
\end{quote}

\begin{center}
\begin{tabular}{||c||c||c||c||c||c||}
  \hline
  $L$ & $Z(D)$ & $x^2$ & $y^2$ & $G$ & $u^2=g_{o}$ \\
  \hline
  $1$ & $\langle t_1 \rangle$ & 1 & 1 & $G_1$ & 1 \\
  $2$ & $\langle t_1 \rangle$ & $t_1$ & $t_1$ & $G_2$ & $t_1$ \\
  $3$ & $\langle t_1 \rangle \times \langle t_2 \rangle$ & 1 & $t_2$ & $G_3$ & 1 \\
  $4$ & $\langle t_1 \rangle \times \langle t_2 \rangle$ & $t_1$ & $t_2$ & $G_4$ & $t_1$ \\
  $5$ & $\langle t_1 \rangle \times \langle u_1 \rangle$ & 1 & $u_1$ & $G_5$ & 1 \\
  $6$ & $\langle t_1 \rangle \times \langle u_1 \rangle$ & $t_1$ & $u_1$ & $G_6$ & $t_1$ \\
  $7$ & $\langle t_1 \rangle \times  \langle t_2 \rangle \times  \langle t_3 \rangle$ & $t_2$ & $t_3$ & $G_7$ & 1 \\
  $8$ & $ \langle t_1 \rangle \times \langle t_2 \rangle \times  \langle t_3 \rangle$ & $t_2$ & $t_3$ & $G_7$ & $t_1$ \\
  $9$ & $\langle t_1 \rangle \times \langle t_2 \rangle \times  \langle t_3 \rangle$ & $t_2$ & $t_3$ & $G_7 \times  \langle t \rangle$ & $t$ \\
  $10$ & $\langle t_1 \rangle \times \langle t_2 \rangle  \times \langle u_1 \rangle$ & $t_2$ & $u_1$ & $G_8$ & 1 \\
  $11$ & $\langle t_1 \rangle \times \langle t_2 \rangle \times \langle u_1 \rangle$ & $t_2$ & $u_1$ & $G_8$ & $t_1$ \\
  $12$ & $\langle t_1 \rangle \times \langle t_2 \rangle \times \langle u_1 \rangle$ & $t_2$ & $u_1$ & $G_8 \times \langle t \rangle$ & $t$ \\
  $13$ & $\langle t_1 \rangle \times \langle t_2 \rangle \times \langle u_1 \rangle$ & $t_2$ & $u_1$ & $G_8 \times \langle w \rangle$ & $w$ \\
  $14$ & $\langle t_1 \rangle \times \langle u_1 \rangle \times \langle u_2 \rangle$ & $u_1$ & $u_2$ & $G_9$ & 1 \\
  $15$ & $\langle t_1 \rangle \times \langle u_1 \rangle \times \langle u_2 \rangle$ & $u_1$ & $u_2$ & $G_9$ & $t_1$ \\
  $16$ & $\langle t_1 \rangle \times  \langle u_1 \rangle \times \langle u_2 \rangle$ & $u_1$ & $u_2$ & $G_9 \times \langle w \rangle$ & $w$ \\
  \hline
  \end{tabular}
\end{center}
\end{theorem}

\emph{Proof:} We will prove, case by case, that each type of
loop listed in those 9 tables of section 3 is listed in the above table. In this proof we will often make use of the fact that if $u,v,w$ are any three elements of $L$ that do not associate then no two distinct elements of $\{u,v,w\}$ commute, $G = \langle u,v,Z(L) \rangle$ is a non abelian group with $G' = L' = \{1,s\}$ and $L = G \cup Gw = M(G,*,w^2)$. A proof of this fact can be founded in [\cite{livro}, Corolary IV.2.3].

Note that the types of loops $L_1$ to $L_4$, $L_7$ to $L_{10}$, $L_{13}$ to $L_{16}$, $L_{19}$ to $L_{22}$ and $L_{37}$ to $L_{40}$ are finite RA loops. So, in [\cite{casofinito}, section 4], it was shown that these loops belongs to the family of loops listed in lines 1,2,3,4,7,8,9 of Table 10, which are the seven families of finite indecomposable RA loops that exist.

Let $L \in L_5 = M(\langle x,y,t_1,w \rangle, *, u^2)$ with $x^2=y^2=1$ and
$u^2=w$. In this case, $ L = M(\langle x,u,t_1,w \rangle, *, y^2)$ is of type 5.

Now, suppose that $L$ is a type of loop listed in $L_6$, i.e.,
$L=M(\langle x,y,t_1,w \rangle,*,u^2)$ with $x^2=y^2=1$, $u^2=t_1w$. Setting
$t_1w = w'$, hence $\langle t_1 \rangle \times \langle w\rangle = \langle t_1\rangle \times \langle w' \rangle$ and in this way, $L=M(\langle x,u,t_1,w'\rangle,*,y^2)$ is of type 5.

If $L=M(\langle x,y,t_1,w \rangle,*,u^2) \in L_{11}$, with $x^2=y^2=t_1$ and
$u^2=w$, we have that $L=M(\langle x,u,t_1,w \rangle,*,y^2)$ is of type 6.

Suppose $L$ is a loop in $L_{12}$. So, $L=M(\langle x,y,t_1,w \rangle,*,u^2)$ where
$x^2=y^2=t_1$ and $u^2=t_1w$. Setting $t_1w=w'$ and changing the
generator $w$ to $w'$, we have that $L=M(\langle x,u,t_1,w' \rangle,*,y^2)$ is of type 6.

Next, if $L=M(\langle x,y,t_1,t_2,w \rangle,*,u^2) \in L_{17}$ where $x^2=1,
y^2=t_2$ and $u^2=w$ then $ L=M(\langle y,u,t_1,t_2,w \rangle,*,x^2)$ is of type 10.

If $L \in L_{18}$, then setting $w' = t_1w$ we have that $L= M(\langle y,u,t_1,t_2,w'\rangle,*,x^2)$ with $y^2 = t_2, u^2 = w', x^2 = 1$ which is a loop of type 10.

Let $L=M(\langle x,y,t_1,t_2,w \rangle,*,u^2) \in L_{23}$,
$x^2=t_1,y^2=t_2,u^2=w$. Then $L=M(\langle y,u,t_1,t_2,w \rangle,*,x^2)$ which is a loop of type 11.

Suppose now that $L$ is loop in $L_{24}$. Setting $t_1w=w'$ we
have that $L=M(\langle y,u,t_1,t_2,w' \rangle,*,x^2)$ is a loop of type 11.

Observe that those loops in $L_{25}$ are those of type 5.

Now, if $L \in L_{26}$ then $L=M(\langle x,y,t_1,u_1 \rangle, *,u^2) =
M(\langle y,u,t_1,u_1 \rangle,*, x^2) \in L_{31}$ and in what follows, we will see that those loops in $L_{31}$ are of type 5 or 6.

If $L \in L_{27}$, then $L=M(\langle x,y,t_1,u_1,t \rangle,*,u^2) =
M(\langle u,y,t_1,u_1,t \rangle,*,x^2)$ which is of type 10.

Let $L$ be a loop in $L_{28}$. If $o(t_1)\leq o(t)$ then
$\langle t_1 \rangle \times \langle u_1 \rangle \times \langle t \rangle = \langle t_1 \rangle \times \langle u_1 \rangle \times
\langle t_1t \rangle$.  Setting $t_1t = t'$, we have
$L=M(\langle x,y,t_1,u_1,t' \rangle,*,u^2) \in L_{27}$ which we had showed before that this loop is of type 10; if $o(t_1) > o(t)$ then $s = t_1^{2^{m_1-1}} =
(t_1t)^{2^{m_1-1}}$ and setting $t_1t = {t_1}'$ we have that
$L=M(\langle x,y,{t_1}',u_1,t \rangle,*,u^2) = M(\langle x,y,{t_1}',u_1 \rangle,*,u^2) \times
\langle t \rangle$  which is not indecomposable.

If $L \in L_{29}$ then $L=M(\langle x,y,t_1,u_1,w \rangle,*,u^2) =
M(\langle u,y,t_1,u_1,w \rangle,*,x^2)$ which is of type 14.

Next, if $L \in L_{30}$, setting $t_1w = w'$, we see that
$L=M(\langle u,y,t_1,u_1,w'\rangle,*,x^2)$ which is of type 14.

Suppose $L \in L_{31}$. If $m_1 = 1$ then changing $x$ to $xu$ we have
that $(xu)^2 = {t_1}^2$ and using Lemma \ref{expoente} we can suppose that
$x^2=1$. Hence $L=M(\langle x,y,t_1,u_1 \rangle ,*,u^2)$ with
$x^2=1,y^2=u_1,u^2=1$ which implies that $L$ is of type 5. If
$m_1 >1$ then changing $u$ by $xu$ we have that $(xu)^2 = {t_1}^{\alpha}$ with $\alpha$ ímpar. So, using Lemma \ref{expoente} we can suppose that $(xu)^2 = t_1$ which implies that
$L$ is of type 6.

Observe that those loops in $L_{32}$ are loops of type 6.

If $L=M(\langle x,y,t_1,u_1,t \rangle,*,u^2) \in L_{33}$ then $L =
M(\langle u,y,t_1,u_1,t \rangle,*,x^2)$ which is of type 11.

Now, suppose that $L=M(\rangle x,y,t_1,u_1,t \rangle,*,u^2) \in L_{34}$ where $x^2=t_1,
y^2=u_1$ and $u^2=t_1t$. If $o(t)\geq o(t_1)$ then we change $t$ to $t'=t_1t$; in this case $L \in L_{33}$ which we had shown that these loops are of type 11. If $o(t) < o(t_1)$ then
changing $x$ to $xu$ we see that $(xu)^2=t$ and setting $t_1t =
{t_1}'$ implies that $L=M(\langle x,y,t_1,u_1,t \rangle,*,u^2)$ with
$x^2=t,y^2=u_1$ and $u^2 = {t_1}'$ which belongs to family
$\mathcal{L}_{11}$.

If $L \in L_{35}$ then $L=M(\langle x,y,t_1,u_1,w \rangle,*,u^2) =
M(\langle y,u,t_1,u_1,w \rangle,*,x^2)$ which is a loop of type 15.

Suppose now that $L=M(\langle x,y,t_1,u_1,w \rangle,*,u^2) \in L_{36}$. Setting
$t_1w = w'$ we have that $L = M(\langle x,y,t_1,u_1,w' \rangle, *, u^2)$ with $x^2=t_1,y^2=u_1$ and $u^2=w'$. Hence
$L_{35}$ are loops of type 15.

If $L \in L_{41}$ then $L=M(\langle x,y,t_1,t_2,t_3,w \rangle,*,u^2) =
M(\langle x,u,t_1,t_2,w,t_3 \rangle,*,y^2)$ which is of type 12.

If $L \in L_{42}$ then changing the generator $w$ to $t_1w = w'$ and
using the same argument as in $L_{41}$ we have that $L_{42}$ is still of type 12.

Note that loops in $L_{43}, L_{44}$ and $L_{45}$ are, respectively, of the types 10,11 and 12.

Let $L=M(\langle x,y,t_1,t_2,u_1,t \rangle,*,u^2) \in L_{46}$. If $o(t) \geq
o(t_1)$ then we can change $t$ to $t_1t$ which implies that $L$ is of type 12. If $o(t) < o(t_1)$ then we can change $t_1$ to $t_1t$ obtaining $L=M(\langle x,y,t_1,t_2,u_1 \rangle,*,u^2) \times \langle t \rangle$ which is not
indecomposable.

Observe that those loops in $L_{47}$ are precisely those of type 13.

If we change the generator $w$ to $ t_1w$ of the loops in
$L_{48}$ we see that these loops are of type 13.

Loops in $L_{49}$ and in $L_{50}$ are, respectively, loops of type 14 and 15.

If $L \in L_{51} = M(\langle x,y, t_1, u_1,u_2, t \rangle, *, u^2)$ with $x^2 = u_1, y^2 = u_2$ and $u^2 = t$ then $L = M(\langle x,u,t_1,u_1,u_2,t \rangle, *, y^2)$ which is a loop of type 13.

Let $L=M(\langle x,y,t_1,u_1,u_2,t \rangle, *,u^2) \in L_{52}$. If $o(t) \geq o(t_1)$
then we can change $t$ to $t_1t$ and
$L= M(\langle u,y,t_1,t,u_1,u_2\rangle,*,x^2)$ which is of type 13. If
$o(t) < o(t_1)$  then changing $t_1$ to $t_1t$ we have that
$L=M(\langle x,y,t_1,u_1,u_2\rangle,*,u^2) \times \langle t \rangle$ which is not indecomposable.

Note that loops in $L_{53}$ are of type 16.

Finally, if we have a loop $L \in L_{54}$, setting $t_1w = w'$, we
see that $L=M(\langle x,y,t_1,u_1,u_2,w' \rangle,*,u^2)$ is a loop of type 16. $\hfill{\fimprova}$

\medskip

We saw that, up to isomorphism, there are at most sixteen types
of indecomposable finitely generated RA loops. Now we will show
that these types of loops listed in the Table 10 are not isomorphic.

\begin{theorem}
The sixteen types of loops listed in the Theorem \ref{classification} are distinct: loops of different types are not isomorphic.
\end{theorem}

\emph{Proof:} Loops of the types 1,2,3,4, 7, 8 and 9 are those of
finite R.A. loops. In \cite{casofinito}, the authors had proved
that they are not isomorphic. Elementary considerations of the
ranks of the centers, shows that we only have to prove that loops of type 4 are not isomorphic to loops of type 6, loops of type 10 are not isomorphic of those in type 11 and loops of type 14 are not isomorphic of those in type 15.

Remember that in $\mathcal{L}_6$ we have that $m_1 =1$. Observe that in $\mathcal{L}_4, u$
is a non central element of order 2. We will show that do no exist
such element in loops of $\mathcal{L}_6$. Suppose that $w =
x^{a}y^{b}t_1^{c}u_1^{d}u^{e} \in L$ where $L$ is a loop of type 6, $ w
\notin Z(L)=Z(G)$ and $w^2=1$. So,
$w^2=x^{2a}y^{2b}u^{2e}s^{ae+be+ab}{t_1}^{2c}{u_1}^{2d} =
{t_1}^{a+e+2c+ae+be+ab}{u_1}^{b+2d} = 1$. Hence $b=d=0$ and
$a+e+ae$ is even which implies $a$ and $e$ are even and then $w \in Z(G)$,
a contradiction.

Now, we are going to do the same with $\mathcal{L}_{14}$ and
$\mathcal{L}_{15}$. In $\mathcal{L}_{14}, u$ is a non central
element of order 2. Let $w =
x^{a}y^{b}{t_1}^{c}{u_1}^{d}{u_2}^{e}u^{f} \in L$ where $L$ is a loop of type 15, $w \notin Z(L)$ such that $w^2=1$. Using that
$x^2=u_1, y^2=u_2$ and $u^2=t_1$ in $\mathcal{L}_{15}$ we have
that $w^2 =
{u_1}^{a+2d}{u_2}^{b+2e}{t_1}^{2^{m_1-1}(ab+af+bf)+f+2c}=1$. Hence
$a=b=d=e=0$ and $f +2c$ is even, which implies $f$ even and then
$w \in Z(G)$, a contradiction.

Similar arguments shows that $\mathcal{L}_{10} \not\simeq
\mathcal{L}_{11}$ and we really have sixteen non isomorphic families of
indecomposable loops. $\hfill{\fimprova}$

\end{document}